\documentclass[11pt]{article}
\usepackage{t1enc}
\usepackage[utf8]{inputenc}
\usepackage{lmodern}
\usepackage{float}
\usepackage{amsmath,amssymb,amsfonts,amsthm,mathrsfs}
\DeclareFontFamily{OT1}{pzc}{}
\DeclareFontShape{OT1}{pzc}{m}{it}{<-> s * [1.2] pzcmi7t}{}
\DeclareMathAlphabet{\mathpzc}{OT1}{pzc}{m}{it}
\usepackage{graphicx}
\usepackage[all,ps]{xy}
\usepackage[colorlinks,urlcolor=cyan,citecolor=blue,linkcolor=blue]{hyperref}

\def\pc{_{\mathrm{prop}}^{\mathrm{comp}}}
\def\Emb{\mathop{\textstyle\rm Emb}}
\def\max{\mathop{\textstyle\rm max}}
\def\prop{\mathrm{prop}}
\def\id{\mathrm{id}}
\newenvironment{prf}%
{\par\noindent\textbf{Proof.\enspace\ignorespaces}}%
{~$\square$\par\medskip}%
{\medskip\par\noindent\textit{Claim.\enspace\ignorespaces}\begin{em}}%
{\end{em}\par}%
{\medskip\par\noindent\textit{Proof.\enspace\ignorespaces}}%
{~$\diamond$\par\medskip}%
{\medskip\par\noindent\textit{Remark.\enspace\ignorespaces}}%
{\par\medskip}%
\newenvironment{ack}%
{\medskip\par\noindent\textbf{Acknowledgement.\enspace\ignorespaces}}%
{\par}%
\newtheorem{thm}{Theorem}[section]%
\newtheorem{claim}[thm]{Claim}%
\newtheorem*{thm*}{Theorem}%
\newtheorem*{claim*}{Claim}%
\theoremstyle{definition}
\newtheorem{defi}[thm]{Definition}%
\newtheorem{rmk}[thm]{Remark}%
\newtheorem{ex}[thm]{Example}%
\title{Stable Pontryagin--Thom construction for \\ proper maps II}
\author{András Csépai}
\date{}

\begin{document}

\maketitle

\let\thefootnote\relax\footnote{\begin{minipage}{0.95\textwidth}
Supported by the ÚNKP-19-2 New National Excellence Program of the Ministry for Innovation and Technology.
~ \includegraphics[height=15pt]{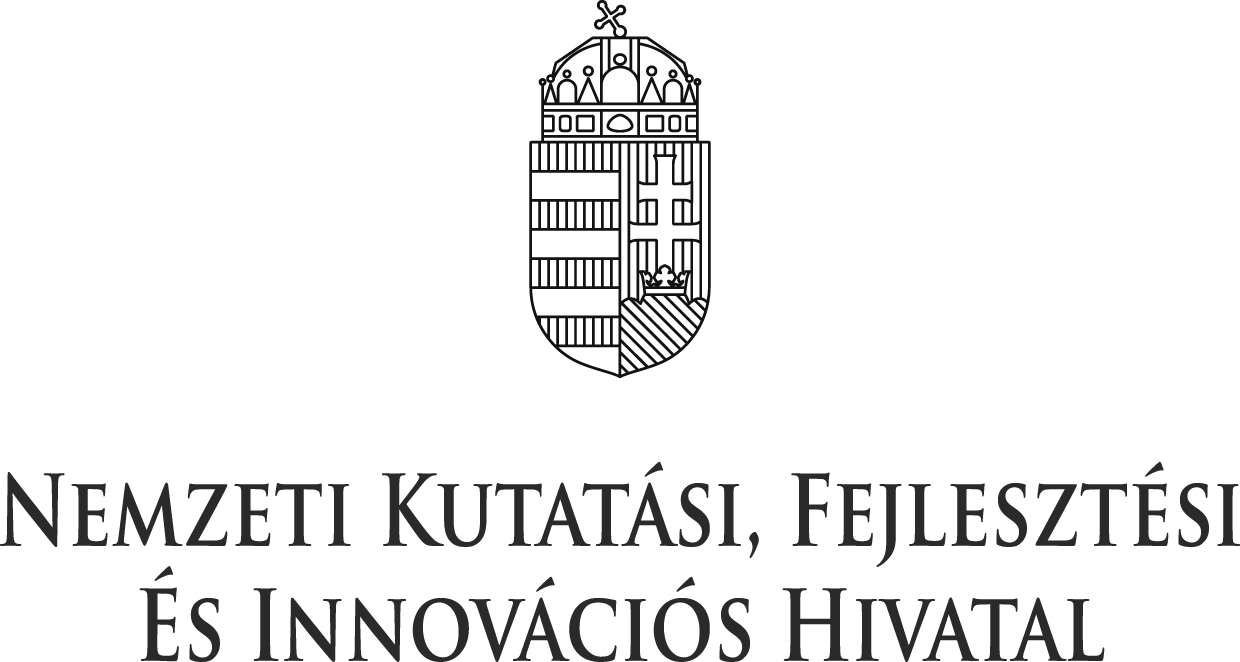} ~ \includegraphics[height=15pt]{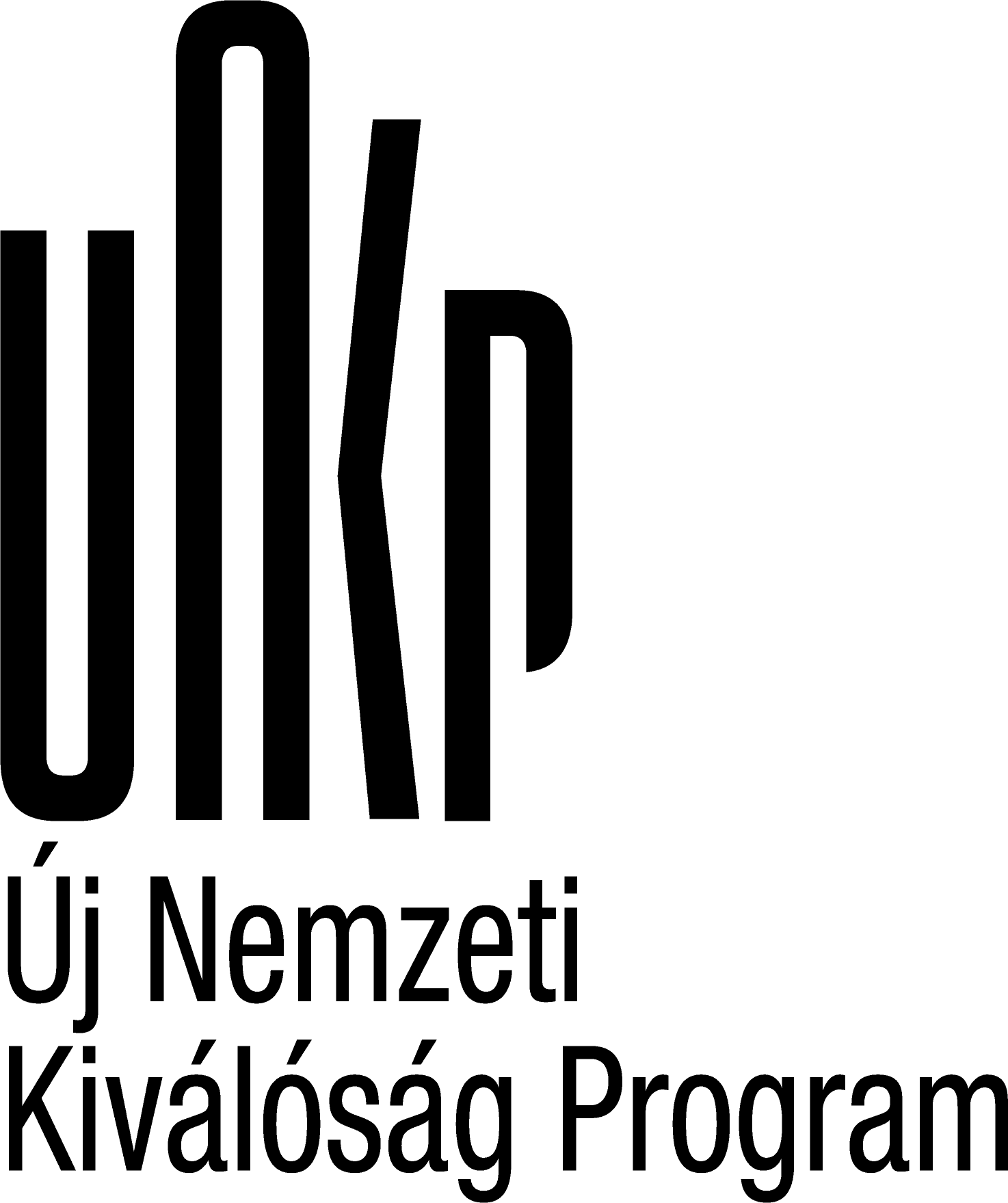}
\end{minipage}}

\begin{abstract}
In \cite{1} we presented a construction that is an analogue of Pontryagin's for proper maps in stable dimensions. This gives a bijection between the cobordism set of framed embedded compact submanifolds in $W\times\mathbb{R}^n$ for a given manifold $W$ and a large enough number $n$, and the homotopy classes of proper maps from $W\times\mathbb{R}^n$ to $\mathbb{R}^{k+n}$. In the present paper we generalise this result in a similar way as Thom's construction generalises Pontryagin's. In other words, we present a bijection between the cobordism set of submanifolds embedded in $W\times\mathbb{R}^n$ with normal bundles induced from a given bundle $\xi\oplus\varepsilon^n$, and the homotopy classes of proper maps from $W\times\mathbb{R}^n$ to a space $U(\xi\oplus\varepsilon^n)$ that depends on the given bundle. An important difference between Thom's construction and ours is that we also consider cobordisms of non-compact manifolds after indroducing a suitable notion of cobordism relation for these.
\end{abstract}

\section{Introduction}

In this paper we will consider cobordisms of submanifolds with a given normal bundle structure in a given manifold and we will establish a connection of these with homotopy classes of so-called proper maps out of the given manifold. This is strongly related to the Pontryagin--Thom construction which we will now recall very briefly.

Pontryagin computed the first two stable homotopy groups of spheres using cobordisms. He did this by constructing an isomorphism of the group $\pi_{m+k}(S^k)$ with the cobordism group of framed embedded $m$-dimensional submanifolds of $S^{m+k}$ (by framed submanifold we mean a submanifold with trivialised normal bundle). This gives rise to the question, what happens if we consider submanifolds in other manifolds and with other types of normal bundles. Namely if we consider the sets $\Emb^\xi(m,W)$ defined as follows.

\begin{defi}
Fix a vector bundle $\xi$ of fibre dimension $k$ and a connected manifold $W$ of dimension $m+k$. Let $M_0$ and $M_1$ be two $m$-dimensional closed embedded submanifolds of $W$ with normal bundles $NM_0$ and $NM_1$ induced from $\xi$. We say that $M_0$ and $M_1$ are $\xi$-cobordant, if there is a compact $(m+1)$-dimensional submanifold with boundary $P\subset W\times[0,1]$ such that $\partial P\subset W\times\{0,1\}$, $\partial P\cap W\times\{0\}=M_0$, $\partial P\cap W\times\{1\}=M_1$ and the normal bundle $NP$ is also induced from $\xi$ and the restriction of $NP$ to the boundary is $NP|_{M_0}=NM_0$ and $NP|_{M_1}=NM_1$. The set of $\xi$-cobordism classes of $m$-dimensional closed submanifolds embedded in $W$ with normal bundles induced from $\xi$ is denoted by $\Emb^\xi(m,W)$.
\end{defi}

Throughout this paper all manifolds and vector bundles are assumed to be smooth.

Thom generalised Pontryagin's construction in the sense that he gave a bijection between $\Emb^\xi(m,W)$ and $[W^*,T\xi]_*$, where the latter denotes the set of based homotopy classes of maps from the one-point compactification of $W$ to the Thom space of $\xi$. It is easy to see that for $W=S^{m+k}$ and $\xi=\varepsilon^k$ (the trivial bundle of dimension $k$) this just gives Pontryagin's bijection.

In the present paper we will consider proper maps and work in the stable case, so let us first describe these.

\begin{defi}
A continuous map $f\colon X\to Y$ is said to be proper if $f^{-1}(C)$ is compact for all compact subsets $C\subset Y$. Two proper maps, $f,g\colon X\to Y$ are called proper homotopic, if there is a proper map $H\colon X\times[0,1]\to Y$ so that $H(\cdot,0)=f$ and $H(\cdot,1)=g$. The proper homotopy classes of proper maps $X\to Y$ will be denoted by $[X,Y]_{\mathrm{prop}}$.
\end{defi}

If $f\colon X\to Y$ is proper, then it is easy to see that the suspension of $f$ defined by
$$Sf\colon X\times\mathbb{R}\to Y\times\mathbb{R};~(x,t)\mapsto(f(x),t)$$
is also proper. Of course this construction can be defined for homotopies as well, so the suspensions of proper homotopic maps are also proper homotopic. Therefore there is a suspension map
$$S\colon[X,Y]_{\mathrm{prop}}\to[X\times\mathbb{R},Y\times\mathbb{R}]_{\mathrm{prop}}.$$

In \cite{1} we presented a construction that gives a bijection for a given manifold $W$ of dimension $m+k$ between the cobordism set $\Emb^{\varepsilon^{k+n}}(m,W\times\mathbb{R}^n)$ and the homotopy set $[W\times\mathbb{R}^n,\mathbb{R}^{k+n}]_{\mathrm{prop}}$ for a sufficiently large $n$. These sets also stabilise as $n\to\infty$ (i.e. if we iterate the suspension map, then it will be bijective after a while) and the ``sufficiently large $n$'' here means that $n$ should be in the stable range. It is proved in \cite{rot} that for $n$ not large enough (that is, not in a stable case) the same bijection is not true. This construction can be thought of as an analogue of Pontryagin's for proper maps. In this paper we generalise it in the same way as Thom generalised Pontryagin's construction.

This paper does not rely on \cite{1} and can be read independently of it.

\begin{ack}
I would like to thank András Szűcs for discussing the topic of this paper with me and for his ideas on cobordisms of open manifolds.
\end{ack}

\section{Preliminaries and formulation of our result}

In our construction we will also consider cobordisms of open (non-compact) manifolds with some restrictions. Since there is no well-known standard definition of cobordism between open manifolds, our first task will be to define such a notion.

\begin{defi}
Fix a (not necessarily compact) submanifold $M$ properly embedded in the connected manifold $W$ (that is, for any compact subset $C\subset W$ the intersection $M\cap C$ is also compact). We say that the normal bundle $NM$ is propely induced from the bundle $\xi$ if the inducing map $M\to B$ to the base space of $\xi$ is a proper map.
\end{defi}

The cobordism of two properly embedded submanifolds with normal bundles properly induced from $\xi$ can be defined in the same way as the cobordism of compact submanifolds. Namely we say that the $m$-dimensional submanifolds $M_0$ and $M_1$ are proper $\xi$-cobordant, if there is a properly embedded $(m+1)$-dimensional submanifold with boundary $P\subset W\times[0,1]$ such that $\partial P\subset W\times\{0,1\}$, $\partial P\cap W\times\{0\}=M_0$, $\partial P\cap W\times\{1\}=M_1$ and the normal bundle $NP$ is also properly induced from $\xi$ and the restriction of $NP$ to the boundary is $NP|_{M_0}=NM_0$ and $NP|_{M_1}=NM_1$. We denote the set of these cobordism classes by $\Emb^\xi(m,W)_\prop$.

However, this definition of cobordism allows too much to change, and in many cases $\Emb^\xi(m,W)_\prop$ will be trivial, as the following easy example shows.

\begin{ex}
Let $\xi$ be the trivial line bundle over $B:=[0,\infty)$ and put $M:=\{0\}\subset\mathbb{R}$ in $W:=\mathbb{R}$. Then clearly $NM$ is properly induced from $\xi$ and the curve $P:=\{(t,\tan(\frac{\pi}{2}t))\mid t\in[0,1)\}$ is properly embedded in $\mathbb{R}\times[0,1]$ and its normal bundle is induced from $\xi$ by a diffeomorphism $P\to[0,\infty)$. Therefore a point in a line is proper null-cobordant.
\end{ex}

This means if we want more interesting cases, then we have to make new restrictions on cobordisms. Then, to compensate these and still get the bijection we are looking for, we also have to make restrictions on homotopies. A good way to do this is to introduce the following ``compact support'' conditions.

\begin{enumerate}
\renewcommand{\labelenumi}{\textrm{\theenumi}}
\renewcommand{\theenumi}{(c\arabic{enumi})}
\item\label{c1} A cobordism $P\subset W\times[0,1]$ between $M_0$ and $M_1$ is up to isotopy such that there is a compact subset $C=C_0\times[0,1]\subset W\times[0,1]$ for which $P\setminus C=(M_0\setminus C_0)\times[0,1]=(M_1\setminus C_0)\times[0,1]$.
\item\label{c2} Assuming the condition \ref{c1}, the inducing map $P\to B$ of the normal bundle $NP$ is such that its restriction to $P\setminus C$ is the fixed homotopy (i.e. the restriction $(M_0\setminus C_0)\times\{t\}\to B$ is the same map for all $t\in[0,1]$).
\item\label{c3} A homotopy $H\colon W\times[0,1]\to X$ is up to isotopy such that there is a compact subset $C_0\subset W$ for which its restriction to $(W\setminus C_0)\times[0,1]$ is the fixed homotopy (i.e. the restriction $(W\setminus C_0)\times\{t\}\to X$ is the same map for all $t\in[0,1]$)
\end{enumerate}

By the term ``up to isotopy'' in \ref{c1} and \ref{c3} we mean that there is a diffeotopy
$$\Phi_t\colon W\times[0,1]\to W\times[0,1]~(t\in[0,1])$$
such that $\Phi_0=\id_{W\times[0,1]}$ and the described condition is true after applying $\Phi_1$. In the case of \ref{c1} this is equivalent to saying that the submanifold $P$ is isotopic through a 1-parameter family of proper embeddings to a submanifold that satisfies the condition. This equivalence is a consequence of an extension to proper embeddings of the usual isotopy extension theorem, which we describe in the appendix.

Now we can define the cobordism and homotopy sets we will use.

\begin{defi}
Fix a vector bundle $\xi$ of fibre dimension $k$ and a connected manifold $W$ of dimension $m+k$. Let $M_0$ and $M_1$ be two $m$-dimensional properly embedded submanifolds of $W$ with normal bundles properly induced from $\xi$. We say that $M_0$ and $M_1$ are compactly supported proper $\xi$-cobordant, if there is a proper $\xi$-cobordism $P\subset W\times[0,1]$ between them that satisfies conditions \ref{c1} and \ref{c2}. The set of compactly supported proper $\xi$-cobordism classes of $m$-dimensional submanifolds properly embedded in $W$ with normal bundles properly induced from $\xi$ is denoted by $\Emb^\xi(m,W)\pc$.
\end{defi}

\begin{defi}
Two proper maps, $f,g\colon W\to X$ are called compactly supported proper homotopic, if there is a proper homotopy $H\colon W\times[0,1]\to X$ between them which satisfies condition \ref{c3}. The compactly supported homotopy classes of proper maps $W\to X$ will be denoted by $[W,X]\pc$.
\end{defi}

In order to make this paper easier to read, we will shorten the terms ``compactly supported proper $\xi$-cobordism'' and ``compactly supported pro- per homotopy'' to the terms ``cobordism'' and ``homotopy'' respectively (the vector bundle $\xi$ will always be clear from context).

Now we define the analogue of the Thom space for proper maps.

\begin{defi}
Let $\xi$ be a vector bundle with base space $B$, total space $E$ and projection $p$. Consider its suspension $\xi\oplus\varepsilon^1$ (which has total space $E\times\mathbb{R}^1$) and define the equivalence relation $\sim$ in the following way: any vector in this bundle is of the form $(b,v,t)$ for $b\in B$, $v\in p^{-1}(b)$ and $t\in\mathbb{R}^1$; let $(b_0,v_0,t_0)\sim(b_1,v_1,t_1)$ iff $v_0=0_{b_0}$, $v_1=0_{b_1}$ and $t_0=t_1\le-1$. The space associated to this bundle is
$$U(\xi\oplus\varepsilon^1):=E\times\mathbb{R}^1/\sim.$$
\end{defi}

For the sake of simplicity we will call the direction of the $\varepsilon^1$ vertical and the direction of $\xi$ horizontal in each fibre of $\xi\oplus\varepsilon^1$ for any bundle $\xi$. The above definition in this terminology means that $U(\xi\oplus\varepsilon^1)$ is the space we get when we identify the ``downwards pointing'' rays $(-\infty,-1]$ of all fibres in the total space of $\xi\oplus\varepsilon^1$.

Our main result is the following.

\begin{thm}\label{t}
For any $m,k\in\mathbb{N}$ there is an $n_0\in\mathbb{N}$ so that for all connected manifolds $W$ of dimension $m+k$, all vector bundles $\xi$ of fibre dimension $k$ and any $n>n_0$, there is a bijection
$$\textstyle\Emb^{\xi\oplus\varepsilon^n}(m,W\times\mathbb{R}^n)\pc\leftrightarrow[W\times\mathbb{R}^n,U(\xi\oplus\varepsilon^n)]\pc.$$
\end{thm}

We remark that the direct analogue of Thom's bijection for proper maps would be a bijection between $\Emb^{\xi}(m,W)\pc$ and $[W,X]\pc$ for a ``proper classifying space'' $X$. However, as we mentioned before, there are counterexamples in \cite{rot} which prove that such a bijection does not hold even in very simple cases. Therefore we apply suspensions to the bundle $\xi$ and to the homotopy classes as well (expecting that these will also stabilise after sufficiently many suspensions) which yields the statement above.


\section{Proof of theorem \ref{t}}

Fix an $n\in\mathbb{N}$, at first without any further condition and a proper map $f\colon W\times\mathbb{R}^n\to U(\xi\oplus\varepsilon^n)$. Note that even though $U(\xi\oplus\varepsilon^n)$ is not necessarily a manifold, it is a manifold around the zero section $B$. Therefore there is a small neighbourhood $V$ of $f^{-1}(B)$ in $W\times\mathbb{R}^n$ such that we can approximate $f$ with a function $g$ that is proper homotopic to $f$, smooth in $V$, transverse to $B$ and $g^{-1}(B)\subset V$. Therefore we may assume that the initial function $f$ had these properties.

Now $M_f:=f^{-1}(B)$ is an $m$-dimensional submanifold of $W\times\mathbb{R}^n$ and it is properly embedded because $f$ is proper. Then the pullback $f^*NB$ of the normal bundle of $B$ in $U(\xi\oplus\varepsilon^n)$ will be the normal bundle $NM_f$ in $W\times\mathbb{R}^n$. But the normal bundle of the zero section in $U(\xi\oplus\varepsilon^n)$ is just $\xi\oplus\varepsilon^n$, therefore $NM_f$ is pulled back from $\xi\oplus\varepsilon^n$. We call $M_f$ the Pontryagin manifold of $f$.

Now we show that the cobordism class of $M_f$ only depends on the homotopy class of $f$ and not the choice of the representative.

\begin{claim}
If $f$ is homotopic to $g$ and $g$ is transverse to the zero section $B$ as well, then $M_f$ is cobordant to $M_g$.
\end{claim}

\begin{prf}
If $H\colon W\times\mathbb{R}^n\times[0,1]\to U(\xi\oplus\varepsilon^n)$ is a homotopy, then we can assume that $H$ is also smooth around the preimage $H^{-1}(B)$ and transverse to $B$. Then it makes sense to talk about the manifold $M_H$ in $W\times\mathbb{R}^n\times[0,1]$, and it is easy to check that $M_H$ is a proper $(\xi\oplus\varepsilon^n)$-cobordism between $M_f$ and $M_g$. The manifold $M_H$ also satisfies conditions \ref{c1} and \ref{c2} because the map $H$ satisfies \ref{c3}, so the cobordism of $M_f$ and $M_g$ is also compactly supported.
\end{prf}

So we have constructed a well-defined map
$$[W\times\mathbb{R}^n,U(\xi\oplus\varepsilon^n)]\pc\to\textstyle\Emb^{\xi\oplus\varepsilon^n}(m,W\times\mathbb{R}^n)\pc.$$
What is left is to construct the inverse for it. In this part of the proof we will need $n$ to be a large number. Later it will be convenient to have $W\times\mathbb{R}^{n+1}$ instead of $W\times\mathbb{R}^n$, so in the remaining part of the proof we will use $n+1$ instead of $n$.

Let $n_0:=\max\{1,m-k+2\}$ so that if $n\ge n_0$, the maps of $(m+1)$-dimensional manifolds into $(m+k+n+1)$-dimensionals can be approximated by embeddings (by Whitney's theorem). Fix an $n\ge n_0$ and an $m$-dimensional properly embedded submanifold $M\subset W\times\mathbb{R}^{n+1}$ such that $NM$ is pulled back from $\xi\oplus\varepsilon^{n+1}$, that is, there is a commutative diagram
\begin{equation}
\begin{aligned}
\xymatrix{
NM\ar[r]^(.43){\tilde{g}}\ar[d] & \xi\oplus\varepsilon^{n+1}\ar[d]\\
M\ar[r]^g & B
}\label{d}
\end{aligned}
\end{equation}
such that $g$ is proper and $\tilde{g}$ is a vector space isomorphism on each fibre. We can of course assume that $g$ and $\tilde{g}$ are smooth. Our aim is to construct a Pontryagin--Thom collapse map $f\colon W\times\mathbb{R}^{n+1}\to U(\xi\oplus\varepsilon^{n+1})$ so that $f$ is proper, the Pontryagin manifold of $f$ is $M$ and $f$ induces the same normal bundle structure as $g$.

Fix a Riemannian metric on $W$ and endow $W\times\mathbb{R}^{n+1}$ with the product metric (where we use the Eucledian metric on $\mathbb{R}^{n+1}$). Now $NM$ is the bundle where the fibre over a point $p\in M$ is the orthogonal complement of $T_pM$ in $T_p(W\times\mathbb{R}^{n+1})$. Using the decomposition $\xi\oplus\varepsilon^{n+1}=\xi\oplus\varepsilon^n\oplus\varepsilon^1$ we can put $N_0M:=g^*(\xi\oplus\varepsilon^n)$, and then we have $NM=N_0M\oplus g^*\varepsilon^1$. Here $g^*\varepsilon^1$ is a trivial line bundle and we can assume that it is pointwise orthogonal to $N_0M$ in $NM$.

We remark that any vector in the bundle $\xi\oplus\varepsilon^{n+1}$ has the form $(b,v,x,t)$ for $b\in B$, $v\in p^{-1}(b)$, $x\in\mathbb{R}^n$ and $t\in\mathbb{R}$, so if we assign to any point $b\in B$ the vector $(b,0_b,0,1)$, we get a trivialisation of the last line bundle $\varepsilon^1$. The pullback of this trivialisation by $g$ will be a normal vector field $u\colon M\to NM$ that trivialises $g^*\varepsilon^1$.

To simplify statements, we will call the last real line in $T_p(W\times\mathbb{R}^{n+1})=T_p(W\times\mathbb{R}^n)\times\mathbb{R}^1$ vertical for all $p\in W\times\mathbb{R}^{n+1}$ and the positive direction on this line will be called upwards.

\begin{claim}\label{cl1}
$M$ is cobordant to a manifold $M'$ such that the normal vector field $u'$ of $M'$ that we get by the same process is vertically upwards in all points of $M'$.
\end{claim}

\begin{prf}
By addendum (i) to the local compression theorem in \cite{rs}, $(M,u)$ can be deformed by an ambient isotopy to a submanifold $(M',u')$, where the normal vector $u'$ is vertically upwards. If we denote the isotopy by
$$\Phi_t\colon W\times\mathbb{R}^{n+1}\to W\times\mathbb{R}^{n+1}~(t\in[0,1])$$
and $M_t:=\Phi_t(M)$, then the normal bundle $NM_t$ is such that $NM=\Phi_t^*NM_t$. The submanifold
$$P:=\bigcup_{t\in[0,1]}M_t\times\{t\}\subset W\times\mathbb{R}^{n+1}\times[0,1]$$
is such that the normal bundle $NP$ is the union of the bundles $NM_t$ for all $t\in[0,1]$. This means that $NP$ is induced from $\xi\oplus\varepsilon^{n+1}$ by the maps $g\circ\Phi_t^{-1}|_{M_t}$, so $P$ is a proper $(\xi\oplus\varepsilon^{n+1})$-cobordism between $M$ and $M'$. To see that the conditions \ref{c1} and \ref{c2} are true for $P$, we just have to apply the diffeotopy of $W\times\mathbb{R}^{n+1}\times[0,1]$ defined for $p\in W\times\mathbb{R}^{n+1}$ and $s\in[0,1]$ by
$$(p,s)\mapsto(\Phi_{(1-s)t}(p),s)~(t\in[0,1]),$$
which finishes the proof.
\end{prf}

Hence we may assume that $u$ was initially vertical.

\begin{claim}\label{cl2}
$M$ is cobordant to a manifold that is in $W\times\mathbb{R}^n\times\{0\}$.
\end{claim}

\begin{prf}
Since $u$ is a vertical normal vector field, the projection of $M$ to $W\times\mathbb{R}^n\times\{0\}$ is an immersion and because of the dimension condition made above we may also assume that it is an embedding.

This also implies that
\begin{alignat*}2
P':=&\{(p_0,tp_1,t)\in W\times\mathbb{R}^{n+1}\times[0,1]\mid\\
&\mid p_0\in W\times\mathbb{R}^n,p_1\in\mathbb{R},(p_0,p_1)\in M,t\in[0,1]\}
\end{alignat*}
is an embedded submanifold of $W\times\mathbb{R}^{n+1}\times[0,1]$. The normal space of $P'$ at each point $(p_0,tp_1,t)$ can be identified with the fibre of $NM$ over $(p_0,p_1)$ by the parallel translation along the line $\mathbb{R}^1$ in $W\times\mathbb{R}^{n+1}\times\{t\}=W\times\mathbb{R}^n\times\mathbb{R}^1\times\{t\}$. Therefore $NP'$ is also induced from $\xi\oplus\varepsilon^{n+1}$, so $P'$ is a proper $(\xi\oplus\varepsilon^{n+1})$-cobordism between $M$ and the projected image of $M$ in $W\times\mathbb{R}^n\times\{0\}$ and conditions \ref{c1} and \ref{c2} can be proved in the same way as in the last proof.
\end{prf}

To summarise the above statements, we can assume that $M$ was initially in $W\times\mathbb{R}^n\times\{0\}$, the normal vector field $u$ is vertically upwards and $N_0M$ is orthogonal to $M$ and $u$ pointwise.

Our plan is now to construct a ``nice'' neighbourhood of $M$, then with the help of this neighbourhood define a map of $W\times\mathbb{R}^{n+1}$ to $U(\xi\oplus\varepsilon^{n+1})$, and then prove that this map satisfies every condition we need.

\begin{claim}\label{tub}
The Riemannian metric on $W$ can be chosen such that there is a small $\varepsilon>0$ so that the exponential map from the normal bundle of $M$ into $W\times\mathbb{R}^{n+1}$ is injective on the $\varepsilon$-neighbourhood of the zero section and the exponential map from $T_p(W\times\mathbb{R}^{n+1})$ into $W\times\mathbb{R}^{n+1}$ is a diffeomorphism on the $\varepsilon$-ball around the origin $0_p$ for all $p\in M$.
\end{claim}

\begin{prf}
We will use that the Riemannian metric on $W$ is complete, and indeed we can assume that we have chosen the metric this way because of the results in \cite{no}. Then all closed and bounded subsets of $W\times\mathbb{R}^{n+1}$ are compact by the Hopf--Rinow theorem. Therefore if we fix a point $*\in W\times\mathbb{R}^{n+1}$, then the closed ball
$$D_r:=\{p\in W\times\mathbb{R}^{n+1}\mid d(*,p)\le r\}$$
is compact for all $r\ge0$, which implies that $M\cap D_r$ is also compact (since $M$ is properly embedded).

Now for any fixed $r\ge0$ there is a small $\varepsilon(r)>0$ so that the exponential map from the normal bundle of $M\cap D_r$ into $W\times\mathbb{R}^{n+1}$ is injective on the $\varepsilon(r)$-neighbourhood of the zero section and the exponential map from $T_p(W\times\mathbb{R}^{n+1})$ into $W\times\mathbb{R}^{n+1}$ is a diffeomorphism on the $\varepsilon(r)$-ball around the origin $0_p$ for all $p\in M\cap D_r$.

We can choose these $\varepsilon(r)$'s such that the function $\varepsilon\colon[0,\infty)\to(0,\infty)$ is smooth and decreasing. If we multiply the metric tensor in all points $p\in\partial D_r$ (i.e. the points for which $d(*,p)=r$) by the positive number $\frac{\varepsilon(0)}{\varepsilon(r)}$ for all $r\ge0$, then we get a new metric tensor which satisfies the properties we need by setting $\varepsilon:=\varepsilon(0)$.
\end{prf}

\begin{rmk}
In the previous proof we have $\frac{\varepsilon(0)}{\varepsilon(r)}\ge1$, so the new distance of any two points is at least their distance in the original metric, hence if the original metric was complete, then so is the new metric. We need this because later we will use the completeness of the Riemannian metric once again.
\end{rmk}

For all $p\in M$ denote by $D_0(p)\subset W\times\mathbb{R}^{n+1}$ the image under the exponential of the $(k+n)$-dimensional open disk of radius $\varepsilon$ around $0_p$ in $N_pM$ orthogonal to $u(p)$ and $T_pM$ (remember that the codimension of $M$ is $k+n+1$, so this disk is well-defined). Define
$$U_0:=\underset{p\in M}{\bigcup}D_0(p),$$
so $U_0$ is a tubular neighbourhood of $M$ in $W\times\mathbb{R}^n\times\{0\}$.

\begin{center}
\begin{figure}[H]
\centering\includegraphics[scale=0.2]{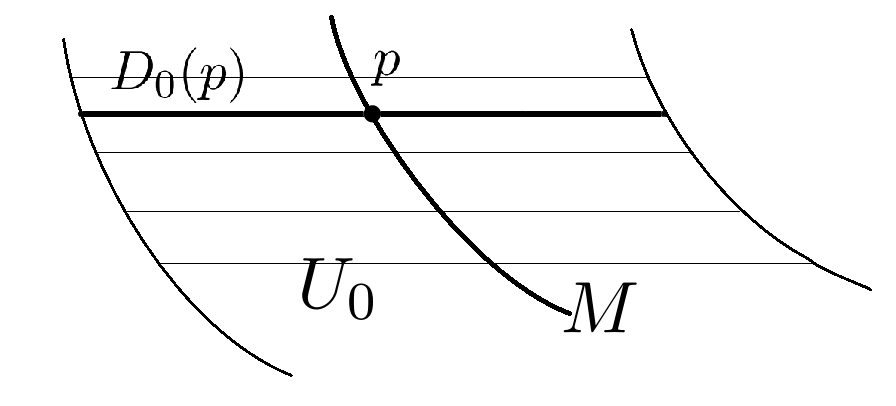}\label{1}
\begin{changemargin}{2cm}{2cm} 
\caption{\hangindent=1.4cm\small We indicated the case when $n=1$ and $k=0$. The arc in the middle represents $M$ and the strip around it is $U_0$. The horizontal segments represent the disks $D_0(p)$ for different $p$'s.}
\end{changemargin} 
\end{figure}
\end{center}
\vspace{-.7cm}

For all $p\in M$ and $q\in D_0(p)$, put $q=(q_0,0)$ where $q_0\in W\times\mathbb{R}^n$. Let $t_q:=-\sqrt{\varepsilon^2-d(p,q)^2}$ (where $d(p,q)$ denotes their distance), so $t_q$ is the negative number for which $d(p,(q_0,t_q))=\varepsilon$. Using the notation $l(q):=\{q_0\}\times(t_q,\infty)$ define
$$U:=\bigcup_{p\in M}\bigcup_{q\in D_0(p)}l(q)=\bigcup_{p\in M}U_p,$$
where $U_p:=\underset{q\in D_0(p)}{\bigcup}l(q)$ is the fibre of $U$ above $p$.

We will also use the following notations: For all points $p\in M$ and $q\in D_0(p)$ we put $l_+(q):=\{q_0\}\times[0,\infty)$ and $l_-(q):=\{q_0\}\times(t_q,0]$. We define
$$D_+(p):=\bigcup_{q\in D_0(p)}l_+(q) ~\text{ and }~ D_-(p):=\bigcup_{q\in D_0(p)}l_-(q),$$
so $D_-(p)$ denotes the half of the $(k+n+1)$-dimensional open disk orthogonal to $M$ with centre $p$ and radius $\varepsilon$, in which the last coordinate of any point is non-positive.

\begin{center}
\begin{figure}[H]
\centering\includegraphics[scale=0.2]{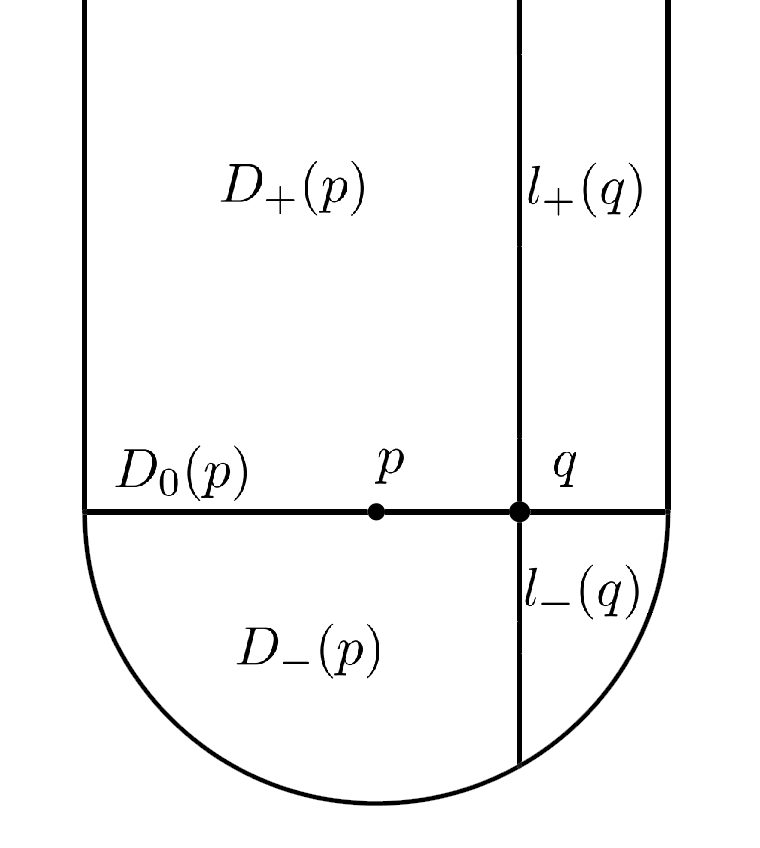}\label{2}
\begin{changemargin}{2cm}{2cm} 
\caption{\hangindent=1.4cm\small $U_p$ looks like this when $n=1$ and $k=0$. We also indicated $l(q)$ for a $q\in D_0(p)$.}
\end{changemargin} 
\end{figure}
\end{center}
\vspace{-.7cm}

\begin{rmk}
It is easy to see that $D_-(p)\cap D_+(p)=D_0(p)$ and $D_-(p)\cup D_+(p)=U_p$ for all $p\in M$ and $U_p$ is diffeomorphic to a neighbourhood of the ray $\{0\}\times[0,\infty)$ in $\mathbb{R}^{k+n+1}$. For different points $p$ the sets $U_p$ are disjoint because $\varepsilon$ was chosen sufficiently small, therefore their union, $U$ is a ``nice'' neighbourhood of $M$.
\end{rmk}

We will also need some new notations in the space $U(\xi\oplus\varepsilon^{n+1})$, so we define these now.

Fix a Riemannian metric on the total space $E$ that is locally the product of a Riemannian metric on $B$ and the Eucledian metric on each fibre of $\xi$. Further, endow the total space $E\times\mathbb{R}^{n+1}$ with the product metric (again we use the Eucledian metric on $\mathbb{R}^{n+1}$).

Take the ($(k+n+1)$-dimensional) disk with centre $(b,0_b,0,-\frac{1}{2})$ and radius $\frac{1}{2}$ in the fibre over $b\in B$ of $\xi\oplus\varepsilon^n\oplus\varepsilon^1$. If we do this for all $b\in B$, then we get a disk bundle which we will denote by $D$ and its boundary sphere bundle which will be denoted by $S$. The fibre of these bundles over $b$ will be $D_b\approx D^{k+n+1}$ and $S_b\approx S^{k+n}$ and the ``upmost'' point of these both is the origin $(b,0_b,0,0)$. 

Let $\downarrow$ denote the downwards unit vector in $U(\xi\oplus\varepsilon^{n+1})$, i.e. $\downarrow$ is the image under the quotient map $E\times\mathbb{R}^n\times\mathbb{R}\to U(\xi\oplus\varepsilon^n\oplus\varepsilon^1)$ of $(b,0_b,0,-1)$ for any $b$. We will denote the images of $D$ and $S$ under this quotient by $D'$ and $S'$ respectively and we remark that these only differ from the bundles $D$ and $S$ in that the ``downmost'' points are identified in all fibres with $\downarrow$.

Now we are ready to define the desired Pontryagin--Thom collapse map for $M$. First fix a point $p\in M$ and define the map on $U_p$.

We map the point $p$ to $g(p)$ in the zero section $B\subset U(\xi\oplus\varepsilon^{n+1})$. Then there is a diffeomorphism
$$f_0\colon D_0(p)\to S_{g(p)}\setminus\{\downarrow\}$$
from the open disk $D_0(p)$ to the punctured sphere $S_{g(p)}\setminus\{\downarrow\}$ that maps the centre $p$ to the north pole $g(p)$. We may choose $f_0$ so that the derivative $df_{0,p}$ maps $T_pD_0(p)$ to $T_{g(p)}S_{g(p)}$ by the same map as $\tilde{g}$ on diagram (\ref{d}). To understand this, notice that $T_pD_0(p)$ is the same as the fibre of $N_0M$ over $p$ and $T_{g(p)}S_{g(p)}$ is the same as the fibre of the subbundle $\xi\oplus\varepsilon^{n}\subset\xi\oplus\varepsilon^n\oplus\varepsilon^1$ over $g(p)$, and by our assumptions on $g$, we have that $\tilde{g}$ maps the fibres of $N_0M$ isomorphically to the fibres of $\xi\oplus\varepsilon^n$.

Then $f_0$ extends to a diffeomorphism
$$f_-\colon D_-(p)\to D_{g(p)}\setminus\{\downarrow\}.$$
This extension can be constructed in the following way: Take a diffeomorphism $\overline{D_-(p)}\approx D^{k+n+1}_-$, where $D^{k+n+1}_-:=\{(x_1,\ldots,x_{k+n+1})\in D^{k+n+1}\mid x_{k+n+1}\le0\}$; compose it with the quotient map
$$D^{k+n+1}_-\to D^{k+n+1}_-/(\partial D^{k+n+1}\cap D^{k+n+1}_-);$$
then identify the quotient space $D^{k+n+1}_-/(\partial D^{k+n+1}\cap D^{k+n+1}_-)\approx D^{k+n+1}$ with $D_{g(p)}$ so that the image of the contracted boundary is $\downarrow$. If we choose these diffeomorphisms so that the restriction of the composed map to $D_0(p)$ is $f_0$, then $f_-$ can be defined as the restriction of this map to $D_-(p)$.

\begin{center}
\begin{figure}[h!]
\centering\includegraphics[scale=0.2]{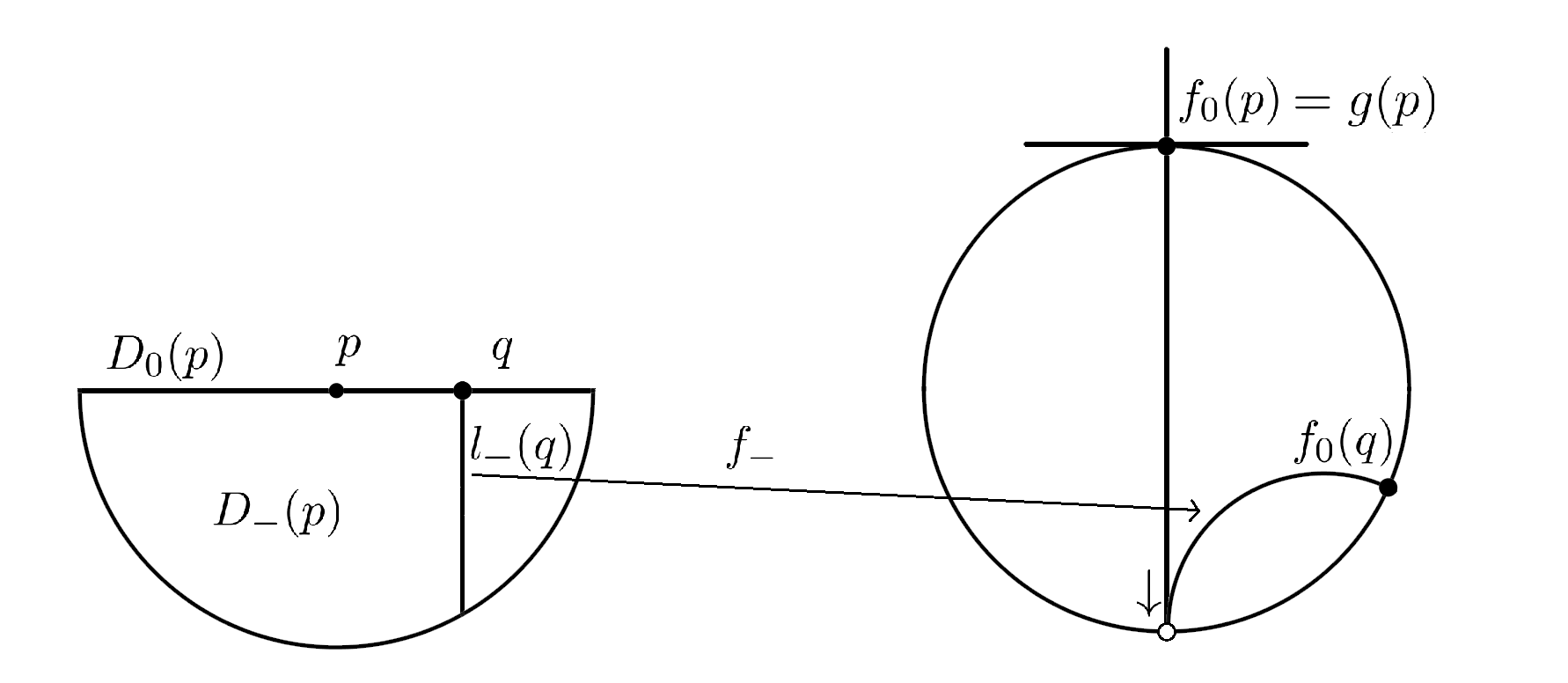}\label{3}
\begin{changemargin}{2cm}{2cm} 
\caption{\hangindent=1.4cm\small The left and the right sides of the figure represent $D_-(p)$ and $D_{g(p)}\setminus\{\downarrow\}$ respectively. The map $f_-$ restricted to $l_-(q)$ is also indicated.}
\end{changemargin} 
\end{figure}
\end{center}
\vspace{-.7cm}

In the following we will use the vectors $\overrightarrow{c,x}$ from the centre $c$ of the sphere $S_{g(p)}$ to any point $x\in S_{g(p)}$. We will also use the rays
$$[x\to):=\{x+t\cdot\overrightarrow{c,x}\mid t\ge0\}$$
and we put $[\downarrow):=[\downarrow\to)$, which is the image under the quotient map $E\times\mathbb{R}^n\times\mathbb{R}\to U(\xi\oplus\varepsilon^n\oplus\varepsilon^1)$ of $[x\to)$ if $x$ is the ``downmost'' point of $S_{g(p)}$.

Take a smooth increasing function $h\colon[0,\infty)\to[0,1]$ such that $h|_{[0,\varepsilon]}\equiv0$ and $h|_{(\varepsilon,\infty)}>0$. For all $q\in D_0(p)$, the ray $l_+(q)$ has a bijection with the ray $[f_0(q)\to)$ defined for $r\in l_+(q)$ by 
$$r\mapsto f_0(q)+(d(r,p)-d(q,p)+h(d(r,p))d(r,*))\cdot\overrightarrow{c,f_0(q)},$$
where $*\in W\times\mathbb{R}^{n+1}$ is a fixed point such that the last real coordinate of $*$ is negative. The union of these maps is a diffeomorphism
$$f_+\colon D_+(p)\to U(\xi\oplus\varepsilon^{n+1})_{g(p)}\setminus((D_{g(p)}\setminus S_{g(p)})\cup[\downarrow))$$
that maps onto the ``fibre'' of $U(\xi\oplus\varepsilon^{n+1})$ minus the open disk $D_{g(p)}\setminus S_{g(p)}$ and the downwards ray $[\downarrow)$. Because of our conditions for $f_-$, the map $f_-\cup f_+$ is a diffeomorphism $U_p\approx U(\xi\oplus\varepsilon^{n+1})_{g(p)}\setminus[\downarrow)$.

\begin{center}
\begin{figure}[h!]
\centering\includegraphics[scale=0.2]{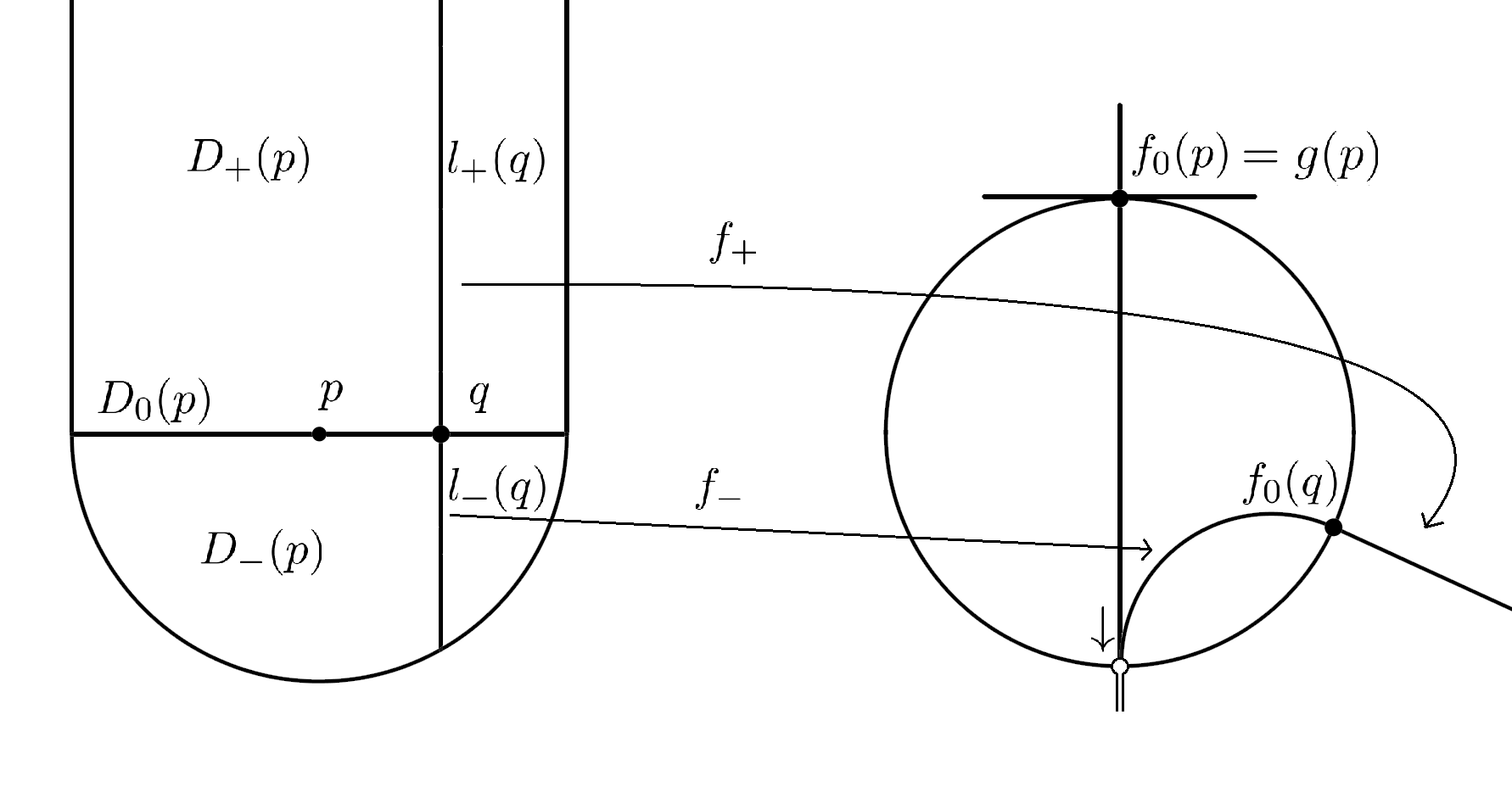}\label{a}
\begin{changemargin}{2cm}{2cm} 
\caption{\hangindent=1.4cm\small The left and the right sides of the figure represent $U_p$ and $U(\xi\oplus\varepsilon^{n+1})_{g(p)}\setminus[\downarrow)$ respectively. The map $f_-\cup f_+$ restricted to $l(q)$ is also indicated.}
\end{changemargin} 
\end{figure}
\end{center}
\vspace{-.7cm}

The above construction defines our Pontryagin--Thom collapse map restricted to any fibre $U_p$. Using the same construction for all $p\in M$, we get the map
$$f_1\colon U\to U(\xi\oplus\varepsilon^{n+1})\setminus[\downarrow).$$
The restriction of $f_1$ to an arbitrary fibre of $U$ is a diffeomorphism, and these diffeomorphisms depend smoothly on the point $p\in M$, because we defined them using the smooth bundle map $\tilde{g}\colon NM\to\xi\oplus\varepsilon^{n+1}$. Hence $f_1$ is smooth.

For all $q\in W\times\mathbb{R}^{n+1}\setminus U$, the distance $d(q,M)$ from $M$ is well-defined and at least $\varepsilon$. Put
\begin{alignat*}2
f_2\colon W\times\mathbb{R}^{n+1}\setminus U & \to[\downarrow)\\
q & \mapsto\downarrow+(d(q,M)-\varepsilon+h(d(q,M))d(q,*))\cdot\overrightarrow{c,\downarrow}.
\end{alignat*}
We remark that this is a well-defined map, i.e. it does not depend on which ``fibre'' we use to take the vector $\overrightarrow{c,\downarrow}$, as this is always just the vector $\overrightarrow{0,-\frac{1}{2}}$ on the vertical line and the rays $[\downarrow)$ are identified in all ``fibres''.

We define the Pontryagin--Thom collapse map as
$$f:=(f_1\cup f_2)\colon W\times\mathbb{R}^{n+1}\to U(\xi\oplus\varepsilon^{n+1}).$$
Now we have to prove that $f$ is a proper map for which $M_f=M$.

\begin{claim}
$f$ is continuous.
\end{claim}

\begin{prf}
First we observe that for all $q\in U_0$ there is a unique $p\in M$ so that $q\in D_0(p)$ and for the point $p$ it holds that $d(q,M)=d(q,p)$. This is because $\varepsilon$ was chosen so that the exponential is a diffeomorphism on the $\varepsilon$-disk for all points of $M$ and $D_0(p)$ is orthogonal to $M$. Then the same is true if $q\in U_p$ is arbitrary, because if $q=(q_0,q_1)$, where $q_0\in W\times\mathbb{R}^n$ and $q_1\in\mathbb{R}$, then
$$d(q,M)=\sqrt{d(q_0,M)^2+q_1^2}=\sqrt{d(q_0,p)^2+q_1^2}=d(q,p).$$
This also implies that if $q\in\overline{U}_p$, then $d(q,M)=d(q,p)$ because the distance is continuous.

It is easy to see that $f_1$ and $f_2$ are both continuous, so we only need to prove that $f$ is continuous in the points of $\partial U$. Choose an arbitrary point $q\in\partial U$ and a sequence $(q_l)$ in $U$ so that $q_l\to q$ as $l\to\infty$. We want to show that the sequence $(f(q_l))$ converges to $f(q)$, or equivalently $(f_1(q_l))$ converges to $f_2(q)$.

There is a $p\in M$ so that $q\in\partial U_p$. If $q\in\partial D_-(p)$, then $d(q,M)=d(q,p)=\varepsilon$, therefore $f_2(q)=\downarrow$. Because of the construction of $f_-$ as a quotient map, $f_1(q_l)\to\downarrow$ as $l\to\infty$, hence $(f_1(q_l))$ indeed converges to $f_2(q)$.

If $q\in\partial D_+(p)$, then we may assume that all of the $q_l$'s are in $\underset{r\in M}{\bigcup}D_+(r)$. The sequence of the unit vectors $\frac{f_1(q_l)}{\lVert f_1(q_l)\rVert}$ converges to $\downarrow$ because $(q_l)$ converges to $\partial U$. By the definition of $f_+$, the norm $\lVert f_1(q_l)\rVert$ tends to
$$\lVert\downarrow\rVert+\lim_{l\to\infty}(d(q_l,M)-d(q_{l,0},M)+h(d(q_l,M))d(q_l,*))\cdot\lVert\overrightarrow{c,\downarrow}\rVert,$$
where $q_{l,0}$ denotes the component of $q_l$ in $W\times\mathbb{R}^n$. The sequence $(q_l)$ converges to $q\in\partial U$, therefore we have $d(q_{l,0},M)\to\varepsilon$, $d(q_l,M)\to d(q,M)$ and $d(q_l,*)\to d(q,*)$ as $l\to\infty$ and $\lVert\downarrow\rVert=1$ and $\lVert\overrightarrow{c,\downarrow}\rVert=\frac{1}{2}$. Hence the sequence of the norms $(\lVert f_1(q_l)\rVert)$ converges to $1+\frac{1}{2}(d(q,M)-\varepsilon+h(d(q,M))d(q,*))$ and so $(f_1(q_l))$ converges again to $f_2(q)$.

Since we have proved this convergence for an arbitrary point and an arbitrary sequence, the continuity of $f$ follows.
\end{prf}

The map $f$ is smooth in a neighbourhood of $M$, $f^{-1}(B)=M$ and $f^*NB=NM=g^*(\xi\oplus\varepsilon^{n+1})$. Therefore if we prove that $f$ is proper, then we get the desired result.

\begin{claim}
$f$ is proper.
\end{claim}

\begin{prf}
Let $C\subset U(\xi\oplus\varepsilon^{n+1})$ be an arbitrary compact subset. Then $C$ is closed and bounded, hence $f^{-1}(C)$ is closed. Put $C=C_1\cup C_2$ where $C_1\subset U(\xi\oplus\varepsilon^{n+1})\setminus[\downarrow)$ and $C_2\subset[\downarrow)$. The set $f_1^{-1}(C_1)$ is bounded, because $C_1$ is bounded and $f_1$ restricted to any fibre of $U$ is the map $f_-\cup f_+$ where taking the preimage $(f_-\cup f_+)^{-1}(C_1)$ only increases the distance of two points by less than $2\varepsilon$, so when we do the same for all points of the compact set $M\cap g^{-1}(C_1)$, we still get a bounded subset (here we used that $g$ is proper). $f_2^{-1}(C_2)$ is trivially bounded, because of the defintion of $f_2$ using the distance from the fixed point $*$ and because $C_2$ is bounded. Hence $f^{-1}(C)$ is a closed and bounded subset of $W\times\mathbb{R}^{n+1}$. If we assume that the Riemannian metric on $W$ is complete, then $f^{-1}(C)$ is compact by the Hopf--Rinow theorem.
\end{prf}

The manifold $M$ was an arbitrary properly embedded submanifold such that an arbitrary proper map $g$ induced $NM$ from $\xi\oplus\varepsilon^{n+1}$, therefore the Pontryagin--Thom construction we have defined assigns to any cobordism class a proper map $f$ for which $M_f$ is in the given cobordism class. Now the only thing left to prove is that the proper homotopy class of $f$ indeed only depends on the cobordism class of $M$ and not the exact representative.

\begin{claim}
If $M$ is cobordant to $M'$, then $f$ is homotopic to the Pontryagin--Thom collapse map $f'$ of $M'$.
\end{claim}

\begin{prf}
If $P$ is a cobordism between $M$ and $M'$, then we can assume that $P$ is compactly supported by conditions \ref{c1} and \ref{c2} (that is, there is a compact set $C\in W\times\mathbb{R}^{n+1}\times[0,1]$ such that $P\setminus C$ is a direct product). By the dimension condition for $n$, all of the constructions made to define $f$ for $M$ have an analogue for $P$, only in the definitions of the maps $f_+$ and $f_2$ we use the map $q\mapsto d(q,\{*\}\times[0,1])$ instead of $q\mapsto d(q,*)$. Therefore there is also a Pontryagin--Thom collapse map for $P$, and it is easy to see that it is a proper homotopy between $f$ and $f'$ which satisfies \ref{c3} since \ref{c1} and \ref{c2} hold for $P$.
\end{prf}

Hence we have an inverse map
$$\textstyle\Emb^{\xi\oplus\varepsilon^{n+1}}(m,W\times\mathbb{R}^{n+1})\pc\to[W\times\mathbb{R}^{n+1},U(\xi\oplus\varepsilon^{n+1})]\pc$$
for $n\ge n_0$, and our proof is complete.

\section{Final remarks}

There are a lot of nice observations concerning theorem \ref{t}, which we would like to collect here.

\begin{rmk}\label{r}
All constructions made above can also be used when we do not assume the compact support conditions (in fact the proof is even easier in this case), which shows that there is a bijection
$$\textstyle\Emb^{\xi\oplus\varepsilon^n}(m,W\times\mathbb{R}^n)_{\prop}\leftrightarrow[W\times\mathbb{R}^n,U(\xi\oplus\varepsilon^n)]_{\prop}.$$
for all manifolds and vector bundles, if $n$ is large enough.
\end{rmk}

\begin{rmk}
In the title of this paper we called our construction ``stable'', which suggests that the sets we use stabilise as $n\to\infty$. We mean by this that the suspension map
$$[W\times\mathbb{R}^n,U(\xi\oplus\varepsilon^n)]\pc\to[W\times\mathbb{R}^{n+1},U(\xi\oplus\varepsilon^{n+1})]\pc$$
is bijective if $n$ is large enough. This is an interesting thing to say, as it is not completely trivial how to define this suspension. As we mentioned in the introduction, suspensions of proper homotopic maps are proper homotopic, but that definition of suspension does not map compactly supported homotopies to compactly supported homotopies.

However, one can define the suspension of cobordism classes by mapping the class represented by $M\subset W$ to the class of $M\times\{0\}\subset W\times\mathbb{R}$ where we add a trivial vertical line bundle to $NM$. Then claims \ref{cl1} and \ref{cl2} imply that this stabilisation is indeed true for the sets $\Emb^{\xi\oplus\varepsilon^n}(m,W\times\mathbb{R}^n)\pc$ (and according to remark \ref{r} also for the sets $\Emb^{\xi\oplus\varepsilon^n}(m,W\times\mathbb{R}^n)_\prop$). Then by theorem \ref{t} this means that the same is true for the sets $[W\times\mathbb{R}^n,U(\xi\oplus\varepsilon^n)]\pc$ (and also for $[W\times\mathbb{R}^n,U(\xi\oplus\varepsilon^n)]_{\mathrm{prop}}$).
\end{rmk}

\begin{rmk}
If the base space $B$ of the bundle $\xi\oplus\varepsilon^n$ is compact, then both $\Emb^{\xi\oplus\varepsilon^n}(m,W\times\mathbb{R}^n)_{\prop}$ and $\Emb^{\xi\oplus\varepsilon^n}(m,W\times\mathbb{R}^n)\pc$ are just the compact cobordism classes of closed submanifolds: $\Emb^{\xi\oplus\varepsilon^n}(m,W\times\mathbb{R}^n)$. Then theorem \ref{t} and remark \ref{r} imply that in this case two maps $W\times\mathbb{R}^n\to U(\xi\oplus\varepsilon^n)$ are compactly supported proper homotopic iff they are proper homotopic.
\end{rmk}

\begin{rmk}
If the bundle $\xi$ is the trivial bundle over a point (that is, $\xi\oplus\varepsilon^n=\varepsilon^{k+n}$), then $U(\varepsilon^{k+n})=\mathbb{R}^{k+n}$. In this case theorem \ref{t} states that $\Emb^{\varepsilon^{k+n}}(m,W\times\mathbb{R}^n)$ is in bijection with $[W\times\mathbb{R}^n,\mathbb{R}^{k+n}]_{\mathrm{prop}}$, which shows that theorem \ref{t} generalises the result proved in \cite{1}.
\end{rmk}

\begin{rmk}
Theorem \ref{t} also implies a nice connection between proper homotopy classes and based homotopy classes between the one-point compactifications. Of course every proper map $X\to Y$ extends to a continuous map $X^*\to Y^*$ by sending infinity to infinity, which defines an injection $[X,Y]_\prop\to[X^*,Y^*]_*$. But the other way is not alwas true, we cannot get all maps $X^*\to Y^*$ as extensions of proper maps $X\to Y$, so there is no inverse to this injection.

However, in our case when the base space $B$ is compact, there are bijections 
\begin{align*}
\xymatrix{
[W\times\mathbb{R}^n,U(\xi\oplus\varepsilon^n)]_{\prop}\ar@{<->}[r] & \Emb^{\xi\oplus\varepsilon^n}(m,W\times\mathbb{R}^n)\ar@{<->}[d]\\
& [(W\times\mathbb{R}^n)^*,T(\xi\oplus\varepsilon^n)]_*
}
\end{align*}
the horizontal one by theorem \ref{t} and the vertical one by the Thom construction. The Thom space of a vector bundle over a compact base space is just the one-point compactification of the total space, so $T(\xi\oplus\varepsilon^n)=(E\times\mathbb{R}^n)^*$. It is easy to see that the one-point compactification of $U(\xi\oplus\varepsilon^n)$ is also homotopy equivalent to $(E\times\mathbb{R}^n)^*$, therefore $[W\times\mathbb{R}^n,U(\xi\oplus\varepsilon^n)]_{\prop}$ is in bijection with $[(W\times\mathbb{R}^n)^*,(U(\xi\oplus\varepsilon^n))^*]_*$.
\end{rmk}

\section*{Appendix: extending isotopies}

The well-known isotopy extension theorem states that an isotopy of a compact embedded submanifold can be extended to a diffeotopy of the ambient manifold. We claim that the same is true if do not assume the submanifold to be compact, only that the isotopy is a 1-parameter family of proper embeddings. The proof we will give is just a slight modification of the proofs in chapter 8 of the textbook \cite{hsch}.

\begin{thm*}
Fix a manifold $W$ and an isotopy $\varphi_t\colon M\to W~(t\in[0,1])$ of the properly embedded submanifold $M\subset W$ such that $\varphi_t$ is a proper embedding for all $t$. Then there is a diffeotopy $\Phi_t\colon W\to W~(t\in[0,1])$ such that $\Phi_t|_M=\varphi_t$ for all $t\in[0,1]$ and $\Phi_0=\id_W$.
\end{thm*}

\begin{prf}
We use claim \ref{tub} to get a complete Riemannian metric on $W$ for which $M\subset W$ has a tubular neighbourhood with radius $\varepsilon>0$. By the compactness of $[0,1]$ we can choose this so that the same $\varepsilon$ works for $\varphi_t(M)$ for all $t\in[0,1]$.

We define the time dependent vector field $X_t\colon W\to TW~(t\in[0,1])$ in the following way: For $p\in\varphi_t(M)$ the vector $X_t(p)$ is the derivative of the curve $s\mapsto\varphi_s(\varphi_t^{-1}(p))$ in $s=t$; for a point $q$ in the tubular neighbourhood of $\varphi_t(M)$ there is a unique $p\in M$ such that $q$ is in the image under the exponential of the $\varepsilon$-disk in $N_p\varphi_t(M)$, the vector $X_t(q)$ is then defined by parallel translating $X_t(p)$ along a minimal geodesic, then multiplying the translated vector by $\frac{\varepsilon-d(q,p)}{\varepsilon}$; outside of the tubular neighbourhood we define $X_t$ to be the zero vector field. We get a continuous vector field this way, but we can approximate it with a smooth vector field so we can assume it was initially smooth.

We will use the maps
$$\tilde{\varphi}_t\colon M\to W\times[0,1];~p\mapsto(\varphi_t(p),t)$$
and define the vector field $Y\colon W\times[0,1]\to T(W\times[0,1])$ by $Y(p,t):=(X_t(p),1)$ for all $t\in[0,1]$. Then all integral curves of $Y$ are parameterised by subintervals of $[0,1]$ and the flow of $Y$ takes the manifold $M=\tilde{\varphi}_0(M)\subset W\times\{0\}$ to $\tilde{\varphi}_1(M)\subset W\times\{1\}$.

\begin{claim*}
All integral curves of $Y$ are compact.
\end{claim*}

If an integral curve is compact, that means it is defined on the closed interval $[0,1]$. If this is true for all integral curves, then they all have the form
$$[0,1]\to W\times[0,1];~t\mapsto(\Phi_t(p),t)$$
for a point $p=\Phi_0(p)\in W$. This defines the map $\Phi_t\colon W\to W$, which is a diffeomorphism for all $t\in[0,1]$ and it is easy to see that $\Phi_t|_M=\varphi_t$ is also true. So the only thing left to prove is the claim above.

\medskip\par\noindent\textbf{Proof of the claim.\enspace\ignorespaces}
If the manifold $M$ is compact, then the vector field $X_t~(t\in[0,1])$ has bounded velocity, which implies that all integral curves of $Y$ have finite length. Then the completeness of the Riemannian metric implies that they are all compact. This also works when $M$ is a manifold with boundary (the vector fields can be constructed in the same way, using the $\varepsilon$-neighbourhood of $M$).

In the general case we take a large compact set $C\subset W$ such that $M\cap C$ is a manifold with boundary. Then if we restrict the isotopy $\varphi_t~(t\in[0,1])$ to $M\cap C$ and construct the vector field $Y'$ for this manifold in the same way as $Y$, then $Y'$ coincides with $Y$ in a neighbourhood of $\varphi_t(M\cap C)$ for all $t$. The compactness of $M\cap C$ implies that the integral curves of $Y'$ are compact and the isotopy of $M\cap C$ extends to a diffeotopy of $W$.

This also implies that the integral curve of $Y'$ starting from a point $q$ in the $\varepsilon$-neighbourhood of $M\cap C$ remains in a small neighbourhood of $\varphi_t(M\cap C)$ for all $t$. If we fix the point $q$ and choose $C$ to be large enough, then the derivative vectors of this integral curve are all vectors of the initial field $Y$. Then the integral curve of $Y$ starting from $q$ is the same as that of $Y'$ because of the uniqueness of the solution for a differential equation.

The reasoning above works for an arbitrary point, so any integral curve of $Y$ is compact. This finishes the proof of the claim and also the proof of the theorem.
\end{prf}

\end{document}